\documentclass[12
pt,twoside]{amsart}
\usepackage{amssymb}
\usepackage{a4wide}
\usepackage[latin1]{inputenc}
\usepackage[T1]{fontenc}
\usepackage{times}

\def\hpq0{h^{p,q}_{\leq 0}}
\def\Hpq0{\H_{\leq 0}^{p,q}}

\def\R{{\mathbb R}}

\def\H{{\mathcal H}}

\def\be{\begin{equation}}
\def\ee{\end{equation}}

\newtheorem{thm}{Theorem}[section]

\newtheorem{cor}[thm]{Corollary}

\theoremstyle{definition}

\theoremstyle{remark}

\newtheorem{preremark}{Remark}
\newtheorem{preex}{Example}

\numberwithin{equation}{section}

\begin{document}

\title[]
{An extension problem for convex functions.}

\author[]{ Bo Berndtsson}

\address{B Berndtsson :Department of Mathematics\\Chalmers University
  of Technology 
  and the University of G\"oteborg\\S-412 96 G\"OTEBORG\\SWEDEN,\\}

\email{ bob@math.chalmers.se}

\begin{abstract}
{ We give a statement on extension with estimates of convex functions defined on a linear subspace, inspired by similar extension results concerning metrics on positive line bundles. }
\end{abstract}

\bigskip

\maketitle

\section{Introduction}
The aim of this paper is to prove the following theorem on extension with estimates of convex functions. 

\begin{thm} Let $\phi(t,x)$ be a convex function in $\R^{n+m}=\R^m_t\times\R^n_x$ and let $\psi(x)$ be a convex function in $\R^n$. Assume
\be
\int_{R^n}e^{\psi(x)-\phi(0,x)} dx =1.
\ee
Then $\psi$ can be extended to a convex function, $\Psi$, on all of $\R^{n+m}$ 
in such a way that
\be
\int_{R^n}e^{\Psi(t,x)-\phi(t,x)} dx \leq 1
\ee
for all $t$ in $\R^m$.
\end{thm}

\bigskip

One motivitation for studying this extension problem comes from the analogy with certain extension problems in complex analytic geometry, in particular the problem of invariance of plurigenera, see e g \cite{Siu} and \cite{Paun}. In these complex analytic extension problems one seeks to extend holomorphic sections to certain line bundles from  hypersurfaces in a bigger manifold.  The main point in the proofs is to extend a positively curved metric on the line bundle, initially defined only over the hypersurface, to a positively curved metric over the ambient manifold. The analogy to the situation discussed in the theorem lies in the parallellism between positively curved metrics and convex functions on $\R^n$. The convex situation is however (and of course) much simpler and allows for more complete  results, and no theorem as precise as theorem 1.1 is known in the complex analytic setting.

\noindent The theorem can also be viewed as a generalization of Prekopa's theorem, \cite{Prekopa},\cite{B-Lieb}, which says that the function $\tilde\phi$ defined by
$$
e^{-\tilde\phi(t)}=\int_{\R^n} e^{-\phi(t,x)} dx
$$
is convex. Indeed, Prekopa's theorem says that if the function $\psi$ is identically equal to zero, then we may take $\Psi(t,x)=\tilde\phi(t)$. 

\noindent It is not hard to see that, conversely, the case $\psi=0$ in theorem 1 implies Prekopa's theorem. This is so because, since $\Psi(t,x)$ is convex with respect to all the variables,
$$
\Psi_0(t):=\inf_x \Psi(t,x)
$$
is also convex, and satisfies (1.2) as well. Hence
$$
\tilde \phi(t)\geq\Psi_0(t)
$$
with equality for $t=0$. In particular, the graph of $\tilde\phi$ has a supporting hyperplane at the origin. Replacing $t=0$ by any other value of $t$, we see that $\tilde\phi$ is convex.

\bigskip

\noindent In the next section we will prove Theorem 1 by a reduction to the case of $\psi=0$. We will also give a simple corollary on the convexity with respect to parameters of a certain extremal convex function.

\bigskip

\noindent
 I would like to thank Mihai P\u aun for many very stimulating discussions on these matters. As mentioned above this little note was inspired by (joint work with him on) similar extension problems for positive metrics on line bundles. Thanks also to the Mittag-Leffler institute where this work was carried out. 

\section{Proof of Theorem 1}

To avoid some issues of convergence we  will first prove a version of 
Theorem 1 where, instead of integrating over $\R^n$, we integrate over a ball $B_R$ in $\R^n$ of radius $R$ and center 0. Let us call a function $\psi$ convex in $B_R$ ``good'', if $\psi$ satisfies the conclusion of Theorem 1 for any choice of $\phi$, convex in $\{|t|<R|\}\times B_R$. More precisely, $\psi$ is good if for any convex function $\phi$ in $\{|t|<R|\}\times B_R$ such that
$$
\int_{B_R} e^{\psi(x)-\phi(0,x)} dx \leq 1,
$$
there is a convex extension $\Psi(t,x)$ such that 

$$
\int_{B_R} e^{\Psi(t,x)-\phi(t,x)} dx \leq 1,
$$
for all $t$  with $|t|<R$.

\bigskip

\noindent By the discussion in the introduction, the function $\psi$ which is identically equal to 0 is good - this is one way of stating Prekopa's theorem. We next claim that any affine $\psi(x)=a\cdot x+b$ is also good. To see this, 
write
$$
1=\int_{B_R} e^{a\cdot x+b-\phi(0,x)} dx =\int_{B_R} 
e^{-(\phi(0,x)-a\cdot x-b)} dx .
$$
Since 0 is good there is a function $\Psi(t,x)$, (or actually $\Psi(t)$), such that $\Psi(0,x)=0$ and
$$
\int_{B_R} e^{\Psi(t,x)-(\phi(t,x)-a\cdot x-b)} dx\leq 1.
$$
Then clearly $\Psi(t,x)+a\cdot x+b$ extends $\psi=a\cdot x+b$ and satisfies the required estimate.

\bigskip

\noindent The next step is to note that if $\psi_\xi(x)$ are good for any $\xi$ in $\R^n$, then $\psi$ defined by
$$
e^\psi =\int e^{\psi_\xi}d\mu(\xi),
$$
where $\mu$ is a positive measure, is also good. This is evident since $\Psi$ defined by
$$
e^{\Psi}=\int e^{\Psi_\xi}d\mu(\xi)
$$
extends $\psi$ if $\Psi_\xi$ extend $\psi_\xi$.

\bigskip

\noindent  The main step of 
the proof involves H\"older's inequality. We claim that if $\psi$ is good and $\lambda\geq 1$, then 
$\psi/\lambda$ is good. This is proved by an iterative procedure, imitating an argument from \cite{B-P}. We can first clearly find an extension $\Psi_0(t,x)$ 
of $\psi$ such that
$$
\int_{B_R} e^{\Psi_0(t,x)/\lambda-\phi(t,x)} dx\leq A,
$$
for {\it some} finite constant $A$. This is at least clear if we shrink $R$ slightly, since we may then take $\Psi_0$ independent of $t$. Write
$$
1=\int_{B_R} e^{\psi(x)/\lambda -\phi(0,x)} dx =\int_{B_R} e^{\psi(x) -
(\phi(0,x)+  
(1-1/\lambda)\Psi_0(0,x))}dx.
$$
Since $\psi$ is good, there is a convex extension $\Psi_1$ with
$$
\int_{B_R} e^{\Psi_1(t,x) -(\phi(t,x)+  (1-1/\lambda)\Psi_0(t,x))} dx\leq 1.
$$
By H\"older's inequality with exponents $\lambda$ and $\lambda/(\lambda-1)$
$$
\int_{B_R} e^{\Psi_1(t,x)/\lambda-\phi(t,x)} dx=
\int_{B_R} e^{\Psi_1(t,x)/\lambda- (1-1/\lambda)\Psi_0(t,x)/\lambda + (1-1/\lambda)\Psi_0(t,x)/\lambda -\phi(t,x)} dx\leq
$$
$$
\leq 
\left(\int_{B_R} e^{\Psi_1(t,x)-(\phi(t,x)+(1-1/\lambda)\Psi_0(t,x))}dx\right)^{1/\lambda} \left(\int_{B_R} e^{\Psi_0(t,x)/\lambda-\phi(t,x)} dx\right)^{(\lambda-1)/\lambda }\leq A^{(\lambda-1)/\lambda}.
$$
If $A>1$ this is strictly smaller than $A$. Iterating the procedure we get extensions $\Psi_k$ of $\psi$ with corresponding integrals bounded by 
$$
A^{((\lambda-1)/\lambda))^k}.
$$
A simple compactness argument shows then that a limit of a subsequence of $\Psi_k$ satisfies the desired estimate. 

\bigskip

\noindent
With this, we can at last prove that any convex function $\psi$ in $\R^n$ is good. Let $\psi^*$ be the Legendre transform
$$
\psi^*(\xi)=\sup_x (x\cdot \xi -\psi(x)).
$$
Then, by the involutivity of the Legendre transform
$$
\psi(x)=\sup_\xi (x\cdot \xi -\psi^*(\xi)).
$$
Hence
$$
\psi=\lim \psi_\lambda/\lambda,
$$
as $\lambda$ tends to infinity, where
$$
e^{\psi_\lambda(x)}=\int_{\R^n}e^{\lambda(x\cdot\xi-\psi^*(\xi))}d\xi.
$$
By the arguments above, each $\psi_\lambda/\lambda$ is good, so by a simple passage to the limit, $\psi$ is good. Finally, we can let $R$ tend to infinity, so the theorem is proved.

\section{ An extremal convex function}

Given a convex function $\phi$ on $\R^n$ we put
$$
E(\phi)(x):=\sup \{\psi(x); \int e^{\psi(x)-\phi(x)}dx \leq 1\}
$$
We then have the following corollary to Theorem 1.1.
\begin{cor} Let $\phi(t,x)$ be convex in $\R^m_t\times \R^n_x$. Let
$$
\hat\phi(t,x)= E_x(\phi)
$$
where $E_x$ indicates that $E$ is taken with respect to the $x$-variable for $t$ fixed. Then $\hat\phi$ is convex.
\end{cor}

\bigskip

Let us show how the Corollary follows from Theorem 1.1. It is  enough
to show that for any point $p=(t_0,x_0)$ there is a convex function $\phi_p$ such 
that
$$
\phi_p\leq\hat \phi,
$$
with equality at $p$. (That implies that the graph of $\hat\phi$ has a supporting hyperplane at every point.) Assume without loss of generality that $p=0$, and let $\psi(x)$ be a function realising the supremum, so that
$$
\psi(0)=\hat\phi(0,0),
$$
and 
$$
\int e^{\psi(x)-\phi(0,x)}dx\leq 1.
$$
By Theorem 1.1 there is a convex function $\Psi(t,x)$ such that $\Psi(0,x)=\psi(x)$ and
$$
\int e^{\Psi(t,x)-\phi(t,x)}dx\leq 1,
$$
for any $t$. Then

$$
\Psi(t,x)\leq \hat\phi(t,x)
$$
for any $(t,x)$ with equality at the origin. This finishes the proof.

\bigskip

\def\listing#1#2#3{{\sc #1}:\ {\it #2}, \ #3.}


\begin{thebibliography}{9999}

 
\bibitem{B-P}\listing{Berndtsson, B and Paun M}{ A Bergman kernel proof of the Kawamata subadjunction theorem.}{arXiv:0804.3884 }



\bibitem{B-Lieb}\listing  {H J Brascamp and E H  Lieb} {On
  extensions of the Brunn-Minkowski and Pr\'ekopa-Leindler theorems,
  including inequalities for log concave functions, and with an
  application to the diffusion equation.}  {J. Functional Analysis  22
  (1976), no. 4, 366--389.}
 
\bibitem{Paun}\listing{ Paun M}{Siu's invariance of plurigenera; a one-tower proof}{ J. Differential Geom.  76  (2007),  no. 3, 485--493}

\bibitem{Prekopa}\listing{Prekopa, A}{On logarithmic concave measures and functions }{ Acad. Sci. Math. (Szeged) 34 (1973), p. 335--343}


\bibitem{Siu}\listing{Siu, Y-T}{ Extension of twisted pluricanonical sections with plurisubharmonic weight and invariance of semipositively twisted plurigenera for manifolds not necessarily of general type}{Complex geometry (G\"ottingen, 2000), 223--277, Springer, Berlin, 2002}

\end{thebibliography}
\end{document}